\documentclass{amsart}

\newtheorem{theorem}{Theorem}
\newtheorem*{theorem*}{Theorem}
\theoremstyle{remark}
\newtheorem{remark}[theorem]{Remark}
\numberwithin{equation}{section}
\usepackage{amssymb}

\def\halfskip{\vskip 10pt plus 1pt minus 1pt}
\def\BB{\mathbb B}
\def\CC{\mathbb C}
\def\GG{\mathbb G}
\def\NN{\mathbb N}
\def\too{\longrightarrow}
\def\CO{\mathcal O}
\def\Aut{\operatorname{Aut}}
\def\const{\operatorname{const}}
\def\id{\operatorname{id}}

\begin{document}

\title
{On automorphisms of the symmetrized bidisc}

\author{Marek Jarnicki}
\address{Jagiellonian University, Institute of Mathematics,
Reymonta 4, 30-059 Krak\'ow, Poland} \email{jarnicki@im.uj.edu.pl}

\author{Peter Pflug}
\address{Carl von Ossietzky Universit\"at Oldenburg, Institut f\"ur Mathematik,
Postfach 2503, D-26111 Oldenburg, Germany}
\email{pflug@mathematik.uni-oldenburg.de}
\thanks{Both authors were supported in part by KBN grant no.~5 P03A 033 21
and by DFG grant no.~227/8-1.}

\subjclass[2000]{32M17}

\keywords{symmetrized bidisc, automorphism}

\begin{abstract}
We characterize the group $\Aut(\GG_2)$ for the symmetrized bidisc
$\GG_2:=\{(\lambda_1+\lambda_2, \lambda_1\lambda_2):
|\lambda_1|, |\lambda_2|<1\}\subset\CC^2$.
\end{abstract}

\maketitle

Let $E$ be the unit disc. Put
\begin{gather*}
\pi:\CC^2\too\CC^2,\quad\pi(\lambda_1,\lambda_2):=
(\lambda_1+\lambda_2, \lambda_1\lambda_2),\\
\GG_2:=\pi(E^2)=\{(\lambda_1+\lambda_2, \lambda_1\lambda_2):
\lambda_1,\lambda_2\in E\},\\
\Sigma_2:=\{\pi(\lambda,\lambda):\lambda\in E\}=
\{(2\lambda, \lambda^2): \lambda\in E\},\\
h_a(\lambda):=\frac{\lambda-a}{1-\overline a\lambda},
\quad a\in E,\;\lambda\in\CC\setminus\{1/\overline a\}.
\end{gather*}
The domain $\GG_2$ is called the {\it symmetrized bidisc}.
It has been recently studied by many authors,
e.g.~\cite{AglYou2001}, \cite{AglYou2003}, \cite{Cos2003a},
\cite{PflZwo2003b}. The domain $\GG_2$ is the first known bounded pseudoconvex
domain for which the Carath\'eodory and Lempert functions coincide, but
which cannot be exhausted by domains biholomorphic to convex domains;
cf.~\cite{Cos2003a}, see also \cite{JarPfl2004} for details.

\begin{remark}
If $f\in\CO(E^2)$ is symmetric, then the relation
$F(\pi(\lambda_1,\lambda_2))=f(\lambda_1,\lambda_2)$ defines a
function $F\in\CO(\GG_2)$.
In particular, if $h\in\CO(E,E)$, then the relation
$H_h(\pi(\lambda_1,\lambda_2))=\pi(h(\lambda_1),h(\lambda_2))$
defines a holomorphic mapping $H_h:\GG_2\too\GG_2$ with
$H_h(\Sigma_2)\subset\Sigma_2$.

If $h\in\Aut(E)$, then  $H_h\in\Aut(\GG_2)$,
$H_h^{-1}=H_{h^{-1}}$, and $H_h(\Sigma_2)=\Sigma_2$.
In particular, if $h(\lambda):=\tau\lambda$ for some
$\tau\in\partial E$, then we get the ``rotation''
$R_\tau(s,p):=H_h(s,p)=(\tau s,\tau^2 p)$.
\end{remark}

\begin{remark}\label{AutactsonS}
For any point $(s_0,p_0)=(2a,a^2)\in\Sigma_2$ we get
$H_{h_a}(s_0,p_0)=(0,0)$. Consequently, the group $\Aut(\GG_2)$
acts transitively on $\Sigma_2$.
\end{remark}

\begin{theorem*}
$$
\Aut(\GG_2)=\{H_h: h\in\Aut(E)\}.
$$
\end{theorem*}

A description of $\Aut(\GG_2)$ is also announced in \cite{AglYou2003}
for a future paper. Since they write ``we shall identify the automorphisms
of $\GG_2$ and this will enable us to show that $\GG_2$ is inhomogeneous'',
it seems that their method will be different from the one in this note.

\begin{proof} First observe that $\Aut(\GG_2)$ does not act
transitively on $\GG_2$.

Otherwise, by the Cartan classification theorem
(cf.~\cite{Akh1990}, \cite{Fuks1965}), $\GG_2$ would be
biholomorphic to the unit Euclidean ball $\BB_2$ or to the unit bidisc $E^2$.
In the case where $\GG_2\simeq \BB_2$ we get a contradiction with the
Remmert--Stein theorem (cf.~\cite{Nar1971}, p.~71) saying that there are no
proper holomorphic mappings $E^2\too\BB_2$. In the case where
$\GG_2\simeq E^2$ we get a contradiction with the characterization of proper
holomorphic mappings $F:E^2\too E^2$ (cf.~\cite{Nar1971}, p.~76), which says that
any such a mapping has the form $F(z_1,z_2)=(F_1(z_1),F_2(z_2))$ (up to a
permutation of the variables). See also \cite{Cos2003a}.

\halfskip

Next observe that $F(\Sigma_2)=\Sigma_2$ for every $F\in\Aut(\GG_2)$.

Indeed, let $V:=\{F(0,0): F\in\Aut(\GG_2)\}$. By W.~Kaup's
theorem, $V$ is a connected complex submanifold of $\GG_2$
(cf.~\cite{Kaup1970}). We already know that $\Sigma_2\subset V$
(Remark \ref{AutactsonS}). Since $\Aut(\GG_2)$ does not act
transitively, we have $V\varsubsetneq\GG_2$. Thus $V=\Sigma_2$.

Take a point $(s_0,p_0)=H_h(0,0)\in\Sigma_2$ with $h\in\Aut(E)$
(Remark \ref{AutactsonS}). Then for every $F\in\Aut(\GG_2)$, we
get $F(s_0,p_0)=(F\circ H_h)(0,0)\in V=\Sigma_2$.

\halfskip

By Remark \ref{AutactsonS}, we only need to show
that every automorphism $F\in\Aut(\GG_2)$ with $F(0,0)=(0,0)$ is
equal to a ``rotation'' $R_\tau$. Fix such an $F=(S,P)$.

First observe that $F|_{\Sigma_2}\in\Aut(\Sigma_2)$. Hence
$F(2\lambda,\lambda^2)=(2\alpha\lambda,\alpha^2\lambda^2)$
for some $\alpha\in\partial E$. Taking $R_{1/\alpha}\circ F$
instead of $F$, we may assume that $\alpha=1$. In particular,
$F'(0,0)\left[\begin{matrix}2 \\ 0\end{matrix}\right]
=\left[\begin{matrix}2 \\ 0\end{matrix}\right]$ and, therefore,
$F'(0,0)=\left[\begin{matrix}1 & b \\ 0 & d\end{matrix}\right]$.

For $\tau\in\partial E$ put $G_\tau:=F^{-1}\circ R_{1/\tau}\circ
F\circ R_\tau\in\Aut(\GG_2)$. Obviously, $G_\tau(0,0)=(0,0)$.
Moreover, $G'_\tau(0,0)=\left[\begin{matrix} 1 &  b(\tau-1) \\ 0 &
1\end{matrix}\right]$. Let $G_\tau^n:\GG_2\too\GG_2$ be the
$n$--th iterate of $G_\tau$.
We have $(G_\tau^n)'(0,0)=\left[\begin{matrix}1 & nb(\tau-1)\\
0 & 1\end{matrix}\right]$. Using the Cauchy inequalities, we get
$$
|nb(\tau-1)|\leq\const,\quad n\in\NN,\;\tau\in\partial E,
$$
which implies that $b=0$, i.e.~$F'(0,0)$ is diagonal.

\halfskip

We have $G'_\tau(0,0)=\id$. Hence, by the Cartan
theorem (cf.~\cite{Nar1971}, p.~66), $G_\tau=\id$. Consequently,
$R_\tau\circ F=F\circ R_\tau$, i.e.
$$
(\tau S(s,p),\tau^2P(s,p))=(S(\tau s,\tau^2p),P(\tau s,\tau^2p)),
\quad (s,p)\in\GG_2,\;\tau\in\partial E.
$$
Hence $F(s,p)=(s,p+Cs^2)$.  Since
$F(2\lambda,\lambda^2)=(2\lambda,\lambda^2)$, we have
$(2\lambda,\lambda^2+4C\lambda^2)=(2\lambda,\lambda^2)$, which
immediately implies that $C=0$, i.e.~$F=\id$.
\end{proof}


\bibliographystyle{amsplain}

\end{document}